\documentclass{amsart}
\usepackage{amsmath}
\usepackage{amsthm}
\usepackage{amssymb}
\usepackage{epsfig}

\theoremstyle{plain}
\newtheorem{thm}{Theorem}[section]
\newtheorem{lem}[thm]{Lemma}
\newtheorem{prop}[thm]{Proposition}
\newtheorem{cor}[thm]{Corollary}
\theoremstyle{definition}
\newtheorem{example}[thm]{Example}
\newtheorem{defn}[thm]{Definition}
\theoremstyle{remark}
\newtheorem{rem}[thm]{Remark}

\newenvironment{myquote}
   {\begin{list}{}{\setlength{\leftmargin}{.35in}\setlength{\rightmargin}{0in}}%
    \item\relax}
   {\end{list}}

\newcommand{\bm}[1]{\mbox{\boldmath $#1$}}

\title{Reduced Decompositions and Permutation Patterns}
\author{Bridget Eileen Tenner}
\address{Department of Mathematics, Massachusetts Institute of Technology, 77 Massachusetts Ave, Cambridge, MA 02139, USA}
\email{bridget@math.mit.edu}
\date{February 17, 2006}

\begin{document}

\begin{abstract}

Billey, Jockusch, and Stanley characterized $321$-avoiding permutations by a property of their reduced decompositions.  This paper generalizes that result with a detailed study of permutations via their reduced decompositions and the notion of pattern containment.  These techniques are used to prove a new characterization of vexillary permutations in terms of their principal dual order ideals in a particular poset.  Additionally, the combined frameworks yield several new results about the commutation classes of a permutation.  In particular, these describe structural aspects of the corresponding graph of the classes and the zonotopal tilings of a polygon defined by Elnitsky that is associated with the permutation.

\end{abstract}

\maketitle

\section{Introduction}\label{section intro}

Reduced decompositions of permutations are classical objects in combinatorics that appear throughout the literature.  Following the work of Rodica Simion and Frank Schmidt in \cite{simion}, the study of permutation patterns, particularly pattern avoidance, has become a frequently studied field as well.

In \cite{bjs}, Sara Billey, William Jockusch, and Richard Stanley relate these two concepts, possibly for the first time.  There they show that $321$-avoiding permutations are exactly those permutations where the subsequence $i (i\pm 1) i$ never occurs in a reduced decomposition.  Relatedly, Victor Reiner shows in \cite{reiner} that the number of $i (i\pm 1)i$ occurrences in reduced decompositions of the longest element in the symmetric group, which has the maximal number of occurrences of $321$, is equal to the number of such reduced decompositions.  Stanley had previously shown that this is the number of standard Young tableaux of a staircase shape in \cite{stanley}.

Inspired by these results, and more generally by the relationship they suggest between the two aspects of permutations, this paper studies elements of the symmetric group from the combined perspectives of their reduced decompositions and their patterns.  While these aspects of a permutation appear extensively in combinatorial literature, they are not often treated together.  This paper strives to remedy that fact, addressing several questions where reduced decompositions and permutation patterns together lead to interesting results.

After introducing basic terminology and notation in Section~\ref{section defs}, Section~\ref{section vexillary} generalizes the result of Billey, Jockusch, and Stanley, via a new characterization of vexillary permutations in Theorem~\ref{vexthm}.  This characterization is based on the reduced decompositions of the permutations \emph{containing} the permutation in question, and is strikingly different from all previous equivalent characterizations.  In addition to requiring that each of the permutations containing the vexillary permutation has a certain kind of reduced decomposition, the proof of Theorem~\ref{vexthm} explicitly constructs such a reduced decomposition.

There are three algorithms which appear in this paper, the first of which occurs in the proof of Theorem~\ref{vexthm}.  It should be noted that these are not deterministic, and include a certain amount of choice.  For instance, Example~\ref{vex ex} describes only one possible route that the algorithm \textsf{VEX} may take on a particular input.

There is an equivalence relation, sometimes known as the commutation relation, on the set of reduced decompositions of a particular permutation.  This and an associated graph are discussed in Section~\ref{section eq}.  Theorem~\ref{1 lbm} and Corollary~\ref{1 ulbm} characterize permutations with graphs and commutation classes having certain properties.  These results are strengthened in Theorem~\ref{chainthm}.

The results in Sections~\ref{section polygon} and \ref{section poset} 
discuss permutation patterns with respect to a polygon defined by Serge 
Elnitsky in \cite{elnitsky}.  The rhombic tilings of this polygon are in 
bijection with the commutation classes of a permutation.  New results 
include that the number of commutation classes of a permutation is 
monotonically increasing with respect to pattern containment 
(Theorem~\ref{monotone}), and several results pertaining to a poset 
associated with tilings of the polygon.  Finally, Section~\ref{section fb} 
completely describes this poset in the case of a freely braided 
permutation, as defined by Richard Green and Jozsef Losonczy in 
\cite{green1} and \cite{green2}.

\section{Basic Definitions}\label{section defs}

The main definitions and notation that appear throughout the paper are discussed below.  For more information about these objects, including proofs of elementary facts, see \cite{bjornerbrenti} and \cite{macdonald}.

Let $\mathfrak{S}_n$ denote the symmetric group on $n$ elements.  An element $w \in \mathfrak{S}_n$ permutes $\{1, \ldots, n\}$ by mapping $i \mapsto w(i)$.  This permutation will be written in one-line notation $w = w(1)w(2) \cdots w(n)$.

\begin{example}
$4213 \in \mathfrak{S}_4$ maps $1$ to $4$, $2$ to itself, $3$ to $1$, and $4$ to $3$.
\end{example}

For $i \in \{1, \ldots, n-1\}$, the map $s_i$ transposes $i$ and $i+1$, and fixes all other elements in a permutation.  The symmetric group $\mathfrak{S}_n$ is the Coxeter group of type $A_{n-1}$, and it is generated by the adjacent transpositions $\{s_i: i = 1, \ldots, n-1\}$.  The adjacent transpositions satisfy the Coxeter relations:
\begin{eqnarray}
& s_i^2 = 1 &\text{ for all } i;\nonumber\\
& s_is_j = s_js_i &\text{ if } |i-j| > 1;\label{shortbraid}\\
& s_is_{i+1}s_i = s_{i+1}s_is_{i+1} &\text{ for } 1 \le i \le n-2.\label{longbraid}
\end{eqnarray}

\noindent Equation~\eqref{shortbraid} is called the \emph{short braid relation}, and equation~\eqref{longbraid} is the \emph{long braid relation}.  A map is written to the left of its input, so $s_iw$ interchanges the positions of values $i$ and $i+1$ in the permutation $w$, while $ws_i$ interchanges the values in positions $i$ and $i+1$ in $w$.  If $w = w(1) \cdots w(n)$, then $ws_i = w(1) \cdots w(i+1) w(i) \cdots w(n)$.

Because the symmetric group is generated by adjacent transpositions, any permutation $w \in \mathfrak{S}_n$ can be written as $w = s_{i_1} \cdots s_{i_\ell}$ for some $\{i_1, \ldots, i_\ell\}$.  The least such $\ell$ is the \emph{length} of $w$, denoted $\ell(w)$.  An \emph{inversion} in $w$ is a pair $(i,j)$ where $i < j$ and $w(i) > w(j)$.  The \emph{inversion set} is $I(w) = \{(i,j): (i,j) \text{ is an inversion}\}$.  Since $I(w) \subseteq [1,n] \times [1,n]$, the inversion set can also be viewed as an array.  The number of inversions in $w$ is equal to $\ell(w)$ (see \cite{macdonald}).  For obvious reasons, the permutation $w_0 := n \cdots 2 1 \in \mathfrak{S}_n$ is called the \emph{longest element} in $\mathfrak{S}_n$.  

\begin{defn}
For a permutation $w$ with $\ell(w) = \ell$, a string $i_1 \cdots i_\ell$ such that $w = s_{i_1} \cdots s_{i_\ell}$ is a \emph{reduced decomposition} of $w$.  (Some sources call this a \emph{reduced word}.)  The set $R(w)$ consists of all reduced decompositions of $w$.
\end{defn}

\begin{defn}
A \emph{factor} is a consecutive substring of a reduced decomposition.
\end{defn}

Similar to the Coxeter relations, a factor $j_1 j_2$ in a reduced decomposition will be called a \emph{short braid move} if $|j_1 - j_2| > 1$, and a factor $j (j \pm 1) j$ will be called a \emph{long braid move}.  The set $R(w)$ has been studied in various contexts, notably by Stanley in \cite{stanley}.  There, Stanley computes $|R(w)|$ for several classes of permutations in terms of the number of standard Young tableaux of certain shapes.  In the case of a vexillary permutation, this is equal to the number of standard Young tableaux of a single shape $\lambda(w)$ (see also Exercise 7.22 of \cite{ec2}).  The definition of vexillary permutations is postponed until Section~\ref{section vexillary}, where they will be discussed in depth.

\begin{defn}
Let $w = w(1) \cdots w(n)$ and $p = p(1) \cdots p(k)$ for $k \le n$.  The permutation $w$ \emph{contains the pattern $p$} if there exist $i_1 < \cdots < i_k$ such that $w(i_1) \cdots w(i_k)$ is in the same relative order as $p(1) \cdots p(k)$.  That is, $w(i_h) < w(i_j)$ if and only if $p(h) < p(j)$.  If $w$ does not contain $p$, then $w$ \emph{avoids} $p$, or is \emph{$p$-avoiding}.
\end{defn}

Suppose that $w$ contains the pattern $p$, with $\{i_1, \ldots, i_k\}$ as defined above.  Then $w(i_1) \cdots w(i_k)$ is an \emph{occurrence} of $p$ in $w$.  The notation $\langle p(j) \rangle$ will denote the value $w(i_j)$.  If $\overline{p} = p(j) p(j+1) \cdots p(j+m)$, then $\langle \overline{p} \rangle = w(i_j) w(i_{j+1}) \cdots w(i_{j+m})$.

\begin{example}
Let $w = 7413625$, $p = 1243$, and $q = 1234$.  Then $1365$ is an occurrence of $p$, with $\langle 1 \rangle = 1$, $\langle 2 \rangle = 3$, $\langle 4 \rangle = 6$, and $\langle 3 \rangle = 5$.  Also, $\langle 2 4 \rangle = 36$.  The permutation $w$ is $q$-avoiding.
\end{example}

\begin{defn}
Let $w$ contain the pattern $p$, and let $\langle p \rangle$ be a particular occurrence of $p$.  If $w(j) \in \langle p \rangle$, then $w(j)$ is a \emph{pattern entry} in $w$.  Otherwise $w(j)$ is a \emph{non-pattern entry}.  If a non-pattern entry lies between two pattern entries in the one-line notation for $w$, then it is \emph{inside} the pattern.  Otherwise it is \emph{outside} the pattern.  ``Inside'' and ``outside'' are only defined for non-pattern entries.
\end{defn}

\begin{defn}\label{obstruction}
Let $\langle p \rangle$ be an occurrence of $p \in \mathfrak{S}_k$ in $w$.  Suppose that $x$ is inside the pattern, that $\langle m \rangle < x < \langle m+1 \rangle$ for some $m \in [1,k-1]$, and that the values $\{\langle m \rangle, x, \langle m+1 \rangle\}$ appear in increasing order in the one-line notation for $w$.  Let $a, b \in \mathbb{N}$ be maximal so that the values
\[\big\{\langle m-a \rangle, \langle m-a+1 \rangle, \ldots, \langle m \rangle, x, \langle m+1 \rangle, \ldots, \langle m+b-1 \rangle, \langle m+b \rangle\big\}\]
\noindent appear in increasing order in the one-line notation for $w$.  The entry $x$ is \emph{obstructed to the left} if a pattern entry smaller than $\langle m-a \rangle$ appears between $\langle m-a \rangle$ and $x$ in $w$.  Likewise, $x$ is \emph{obstructed to the right} if a pattern entry larger than $\langle m+b \rangle$ appears between $x$ and $\langle m+b \rangle$ in $w$.
\end{defn}

\begin{example}
Let $w = 32451$ and $p = 3241$.  Then $3241$ and $3251$ are both occurrences of $p$ in $w$.  Obstruction is only defined for the latter, with $x=4$ and $m=3$.  Then $a=b=0$, and $4$ is obstructed to the left and not to the right.
\end{example}

\begin{example}\label{nonvex obst}
Let $w = 21354$ and $p = 2143$.  Then $2154$ is an occurrence of $p$ in $w$.  Using $x=3$, $m=2$ in Definition~\ref{obstruction} shows that $a = b = 0$, and $3$ is obstructed both to the left and to the right.
\end{example}

\section{Vexillary Characterization}\label{section vexillary}

Vexillary permutations first appeared in \cite{lascoux} and subsequent publications by Alain Lascoux and Marcel-Paul Sch\"{u}tzenberger.  They were also independently found by Stanley in \cite{stanley}.  There have since emerged several equivalent definitions of these permutations, and a thorough discussion of these occurs in \cite{macdonald}.  The original definition of Lascoux and Sch\"{u}tzenberger, and the one of most relevance to this discussion, is the following.

\begin{defn}\label{2143vex}
A permutation is \emph{vexillary} if it is $2143$-avoiding.
\end{defn}

\begin{example}
The permutation $3641572$ is vexillary, but $3641752$ is not vexillary because $3175$ is an occurrence of $2143$ in the latter.
\end{example}

The following proposition is key to proving one direction of Theorem~\ref{vexthm}.

\begin{prop}\label{vex not obst}
Let $w$ contain the pattern $p$.  Let $x$ be inside the pattern, with $\langle m \rangle < x < \langle m+1 \rangle$ and the values $\{\langle m \rangle, x, \langle m+1 \rangle\}$ appearing in increasing order in $w$.  If $p$ is vexillary then $x$ cannot be obstructed both to the left and to the right.
\end{prop}

\begin{proof}
Such obstructions would create a $2143$-pattern in $p$.
\end{proof}

Example~\ref{nonvex obst} illustrates a non-vexillary permutation which has an element $x$ that is obstructed on both sides.

Equivalent characterizations of vexillarity concern the inversion set $I(w)$ or the following objects.

\begin{defn}
The \emph{diagram} of a permutation $w$ is $D(w) \subseteq [1,n] \times [1,n]$ where
\begin{equation*}
(i,j) \in D(w) \text{ if and only if } i < w^{-1}(j) \text{ and } j < w(i).
\end{equation*}
\end{defn}

\begin{defn}
The \emph{code} of $w$ is the vector $c(w) = (c_1(w), \ldots, c_n(w))$ where $c_i(w)$ is the number of elements in row $i$ of $I(w)$.  The shape $\lambda(w)$ is the partition formed by writing the entries of the code in non-increasing order.
\end{defn}

\begin{prop}\label{vex defs}
The following are equivalent definitions of vexillarity for a permutation $w$:
\begin{enumerate}
\renewcommand{\labelenumi}{(V\arabic{enumi})}
\item\label{vex1} $w$ is $2143$-avoiding;
\item\label{vex2} The set of rows of $I(w)$ is totally ordered by inclusion;
\item\label{vex3} The set of columns of $I(w)$ is totally ordered by inclusion;
\item\label{vex4} The set of rows of $D(w)$ is totally ordered by inclusion;
\item\label{vex5} The set of columns of $D(w)$ is totally ordered by inclusion;
\item\label{vex6} $\lambda(w)' = \lambda(w^{-1})$, where $\lambda(w)'$ is the transpose of $\lambda(w)$.
\end{enumerate}
\end{prop}

\begin{proof}
See \cite{macdonald}.
\end{proof}

This section proves a new characterization of vexillary permutations, quite different from those in Proposition~\ref{vex defs}.  A partial ordering can be placed on the set of all permutations $\mathfrak{S}_1 \cup \mathfrak{S}_2 \cup \mathfrak{S}_3 \cup \cdots$, where $u < v$ if $v$ contains the pattern $u$.  Definition~\ref{2143vex} determines vexillarity by a condition on the principal order ideal of a permutation.  The new characterization, Theorem~\ref{vexthm}, depends on a particular condition holding for the principal \emph{dual} order ideal.

\begin{defn}
Let $\bm{i} = i_1 \cdots i_\ell$ be a reduced decomposition of $w = w(1) \cdots w(n)$.  For $M \in \mathbb{N}$, the \emph{shift of $\bm{i}$ by $M$} is
\[\bm{i}^M := (i_1 + M) \cdots (i_\ell + M) \in R\big(12 \cdots M (w(1)+M)(w(2)+M)\cdots (w(n)+M)\big).\]
\end{defn}

\begin{thm}\label{vexthm}
The permutation $p$ is vexillary if and only if, for every permutation $w$ containing a $p$-pattern, there exists a reduced decomposition $\bm{j} \in R(w)$ containing some shift of an element $\bm{i} \in R(p)$ as a factor.
\end{thm}

\begin{proof}
First suppose that $p \in \mathfrak{S}_k$ is vexillary.  Let $w \in \mathfrak{S}_n$ contain a $p$-pattern.  
Assume for the moment that there is a
\begin{equation}\label{tilde}
\widetilde{w} = \left(s_{I_1}\cdots s_{I_q}\right)w\left(s_{J_1} \cdots s_{J_r}\right) \in \mathfrak{S}_n
\end{equation}
\noindent such that
\begin{enumerate}
\renewcommand{\labelenumi}{(R\arabic{enumi})}
\item $\ell(\widetilde{w}) = \ell(w) - (q + r)$;
\item $\widetilde{w}$ has a $p$-pattern in positions $\{1+M, \ldots, k+M\}$ for some $M \in [0,n-k]$.
\end{enumerate}
Choose a reduced decomposition $\bm{i} \in R(p)$.  Let $\widetilde{w}' \in \mathfrak{S}_n$ be the permutation obtained from $\widetilde{w}$ by placing the values $\{\widetilde{w}(1+M), \ldots, \widetilde{w}(k+M)\}$ in increasing order and leaving all other entries unchanged.  Choose any $\bm{h} \in R(\widetilde{w}')$.  Then
\begin{equation}\label{redw}
\left(I_q \cdots I_1\right) \bm{h} \bm{i}^M \left(J_r \cdots J_1 \right) \in R(w).
\end{equation}
It remains only to find a $\widetilde{w} \in \mathfrak{S}_n$ satisfying (R$1$) and (R$2$).  This will be done by an algorithm \textsf{VEX} that takes as input a permutation $w \in \mathfrak{S}_n$ containing a $p$-pattern and outputs the desired permutation $\widetilde{w} \in \mathfrak{S}_n$.  Because the details of this algorithm can be cumbersome, a brief description precedes each of the major steps.

\medskip

\begin{myquote}

\setlength{\leftmargini}{0in}
\setlength{\leftmarginii}{.18in}
\setlength{\leftmarginiii}{.18in}
\setlength{\leftmarginiv}{.18in}

\noindent \textsf{\textbf{Algorithm} VEX}

\nopagebreak[4]

\noindent \textsf{INPUT: $w \in \mathfrak{S}_n$ with an occurrence $\langle p \rangle$ of the pattern $p \in \mathfrak{S}_k$.}

\noindent \textsf{OUTPUT: $\widetilde{w} \in \mathfrak{S}_n$ as in equation~\eqref{tilde} satisfying (R1) and (R2).}

\smallskip

\begin{enumerate}
\setcounter{enumi}{-1}
\renewcommand{\labelenumi}{\textsf{Step \arabic{enumi}}:}
\renewcommand{\labelenumii}{\textsf{\alph{enumii}}.}
\renewcommand{\labelenumiii}{\textsf{\roman{enumiii}}.}

\item \textsf{Initialize variables.}\\
Set $w_{[0]} := w$ and $i := 0$.
\item\label{alg halt} \textsf{Check if ready to output.}\\
If $w_{[i]}$ has no entries inside the pattern, then \textsf{OUTPUT $w_{[i]}$}.  Otherwise, choose $x_{[i]}$ inside the pattern.  
\item\label{alg big} \textsf{Move all inside entries larger than $\langle k \rangle$ to the right of $\langle p \rangle$.}\\
If $x_{[i]} > \langle k \rangle$, then \textsf{BEGIN}

\begin{enumerate}
\item Let $B(x_{[i]})= \{y \ge x_{[i]} : y \text{ is inside the pattern}\}$.
\item Consider the elements of $B(x_{[i]})$ in decreasing order.  Multiply $w_{[i]}$ on the right by adjacent transpositions (changing \emph{positions} in the one-line notation) to move each element immediately to the right of $\langle p \rangle$.
\item Let $w_{[i+1]}$ be the resulting permutation.  Set $i:=i+1$ and \textsf{GOTO Step~\ref{alg halt}}.
\end{enumerate}

\item\label{alg small} \textsf{Move all inside entries smaller than $\langle 1 \rangle$ to the left of $\langle p \rangle$.}\\
If $x_{[i]} < \langle 1 \rangle$, then \textsf{BEGIN}

\begin{enumerate}
\item Let $S(x_{[i]}) = \{y \le x_{[i]} : y \text{ is inside the pattern}\}$.
\item Consider the elements of $S(x_{[i]})$ in increasing order.  Multiply $w_{[i]}$ on the right by adjacent transpositions to move each element immediately to the left of $\langle p \rangle$.
\item Let $w_{[i+1]}$ be the resulting permutation.  Set $i:=i+1$ and \textsf{GOTO Step~\ref{alg halt}}.
\end{enumerate}

\item\label{alg bound} \textsf{Determine bounds in the pattern for the inside entry.}\\
Let $m \in [1,k-1]$ be the unique value such that $\langle m \rangle < x_{[i]} < \langle m+1 \rangle$.

\item\label{alg between} \textsf{Change the occurrence of $\langle p \rangle$ and the inside entry so that it does not lie between its bounds in the pattern.}\\
If the values $\{\langle m \rangle, x_{[i]}, \langle m+1 \rangle\}$ appear in increasing order in $w_{[i]}$, then define $a$ and $b$ as in Definition~\ref{obstruction} and \textsf{BEGIN}

\begin{enumerate}
\item\label{alg unob right} If $x_{[i]}$ is unobstructed to the right, then \textsf{BEGIN}

\begin{enumerate}
\item Let $R(x_{[i]})$ be the set of non-pattern entries at least as large as $x_{[i]}$ and lying between $x_{[i]}$ and $\langle m+b \rangle$ in the one-line notation of $w_{[i]}$.
\item Consider the elements of $R(x_{[i]})$ in decreasing order.  For each $y \in R(x_{[i]})$, multiply on the right by adjacent transpositions until $y$ is immediately to the right of $\langle m+b \rangle$, or the right neighbor of $y$ is $z > y$.  In the latter case, $z = \langle m+b' \rangle$ for some $b' \in [1,b]$ because all larger non-pattern entries are already to the right of $\langle m+b \rangle$.  Interchange the roles of $y$ and $\langle m+b' \rangle$, and move this new $y$ to the right in the same manner, until it is to the right of (the redefined) $\langle m+b \rangle$.
\item Let $w_{[i+1]}$ be the resulting permutation, with $\langle p \rangle$ 
redefined as indicated.  Let $x_{[i+1]}$ be the non-pattern entry in the final 
move after any interchange of roles.  This is greater than $x_{[i]}$ and the 
newly redefined $\langle m+b \rangle$, and occurs to the right of the new 
$\langle m+b \rangle$.  If $x_{[i+1]}$ is outside of the pattern, \textsf{GOTO 
Step~\ref{alg halt}} with $i:=i+1$.  Otherwise \textsf{GOTO Step~\ref{alg big}} with $i:=i+1$. \end{enumerate}

\item \label{alg unob left} The entry $x_{[i]}$ is unobstructed to the left (Proposition~\ref{vex not obst}).  \textsf{BEGIN}

\begin{enumerate}
\item Let $L(x_{[i]})$ be the set of non-pattern entries at most as large as $x_{[i]}$ and lying between $\langle m-a \rangle$ and $x_{[i]}$ in the one-line notation of $w_{[i]}$.
\item Consider the elements of $L(x_{[i]})$ in increasing order.  For each $y \in L(x_{[i]})$, multiply on the right by adjacent transpositions until $y$ is immediately to the left of $\langle m-a \rangle$, or the left neighbor of $y$ is $z < y$.  In the latter case, $z = \langle m-a' \rangle$ for some $a' \in [0,a]$ because all smaller non-pattern entries are already to the left of $\langle m-a \rangle$.  Interchange the roles of $y$ and $\langle m-a' \rangle$, and move this new $y$ to the left in the same manner, until it is to the left of (the redefined) $\langle m-a \rangle$.
\item Let $w_{[i+1]}$ be the resulting permutation, with $\langle p \rangle$ 
redefined as indicated.  Let $x_{[i+1]}$ be the non-pattern entry in the final 
move after any interchange of roles.  This is less than $x_{[i]}$ and the 
newly redefined $\langle m-a \rangle$, and occurs to the left of the new 
$\langle m-a \rangle$.  If $x_{[i+1]}$ is outside of the pattern, \textsf{GOTO 
Step~\ref{alg halt}}.  Otherwise \textsf{GOTO Step~\ref{alg small}} with $i:=i+1$. \end{enumerate}

\end{enumerate}

\item\label{alg big left} \textsf{Change the occurrence of $\langle p \rangle$, but not its position, so that the value of the inside entry increases but the values of $\langle p \rangle$ either stay the same or decrease.}\\
If $w_{[i]}(s) = \langle m+1 \rangle$ and $w_{[i]}(t) = x_{[i]}$ with $s<t$, multiply $w_{[i]}$ on the left by adjacent transpositions (changing \emph{values} in the one-line notation) to obtain $w_{[i+1]}$ with the values $[x_{[i]}, \langle m+1 \rangle]$ in increasing order.  Then $w_{[i+1]}(s)$ is in the half-open interval $[x_{[i]}, \langle m+1 \rangle)$, and $w_{[i+1]}(t)$ is in the half-open interval $(x_{[i]}, \langle m+1 \rangle]$.  \textsf{GOTO Step~\ref{alg big}} with $x_{[i+1]}:=w_{[i+1]}(t)$, the pattern redefined so that $\langle m+1 \rangle:=w_{[i+1]}(s)$, and $i:=i+1$.
\item\label{alg small right} \textsf{Change the occurrence of $\langle p \rangle$, but not its position, so that the value of the inside entry decreases but the values of $\langle p \rangle$ either stay the same or increase.}\\
If $w_{[i]}(s) = \langle m \rangle$ and $w_{[i]}(t) = x_{[i]}$ with $s>t$, multiply $w_{[i]}$ on the left by adjacent transpositions to obtain $w_{[i+1]}$ with the values $[\langle m \rangle, x_{[i]}]$ in increasing order.  Then $w_{[i+1]}(s)$ is in the half-open interval $(\langle m \rangle, x_{[i]}]$, and $w_{[i+1]}(t)$ is in the half-open interval $[\langle m \rangle, x_{[i]})$.  \textsf{GOTO Step~\ref{alg small}} with $x_{[i+1]}:=w_{[i+1]}(t)$, the pattern redefined so that $\langle m \rangle:=w_{[i+1]}(s)$, and $i:=i+1$.
\end{enumerate}

\end{myquote}

\medskip

Each subsequent visit to Step~\ref{alg halt} involves a permutation with strictly fewer entries inside the pattern than on the previous visit.  Each multiplication by an adjacent transposition indicated in the algorithm removes an inversion, and so decreases the length of the permutation.  This is crucial because of requirement (R$1$). 

Consider the progression of \textsf{VEX}:
\begin{itemize}
\item Step~\ref{alg halt} $\Longrightarrow$ \textsf{HALT} or begin a pass through \textsf{VEX};
\item Step~\ref{alg big} $\Longrightarrow$ Step~\ref{alg halt};
\item Step~\ref{alg small} $\Longrightarrow$ Step~\ref{alg halt};
\item Step~\ref{alg unob right} $\Longrightarrow$ Steps~\ref{alg big} or~\ref{alg big left};
\item Step~\ref{alg unob left} $\Longrightarrow$ Steps~\ref{alg small} or~\ref{alg small right};
\item Step~\ref{alg big left} $\Longrightarrow$ Steps~\ref{alg big}, \ref{alg between}, or~\ref{alg big left};
\item Step~\ref{alg small right} $\Longrightarrow$ Steps~\ref{alg small}, \ref{alg between}, or~\ref{alg small right}.
\end{itemize}

Step~\ref{alg unob right} concludes with $x_{[i+1]}$ to the left of its lower pattern bound, and smaller pattern elements lying between $x_{[i+1]}$ and this bound.  Therefore, no matter how often Step~\ref{alg big left} is next called, the algorithm will never subsequently go to Step~\ref{alg unob left} before going to Step~\ref{alg halt}.  Likewise, a visit to Step~\ref{alg unob left} means that Step~\ref{alg unob right} can never be visited until Step~\ref{alg halt} is visited and a new entry inside the pattern is chosen.

Steps~\ref{alg big} and~\ref{alg small} do not change the relative positions of $\langle p \rangle$.

Steps~\ref{alg unob right} and~\ref{alg big left} imply $x_{[i+1]} > x_{[i]}$, while $x_{[i+1]} < x_{[i]}$ after Steps~\ref{alg unob left} and~\ref{alg small right}.  Let $m$ be as in Step~\ref{alg bound}.  Until revisiting Step~\ref{alg halt}, the values $\langle m' \rangle$, for $m' \ge m+1$, do not increase if  $x_{[i+1]} > x_{[i]}$.  Nor do the values $\langle m' \rangle$, for $m' \le m$, decrease if $x_{[i+1]} < x_{[i]}$.  The other pattern values are unchanged.  The definition of $m$ means that the reordering of values in Steps~\ref{alg big left} and~\ref{alg small right} does not change the positions in which the pattern $p$ occurs.  Additionally, these steps change the value of the entry inside the pattern (that is, $x_{[i+1]} \neq x_{[i]}$), but not its position.

These observations indicate not only that \textsf{VEX} terminates, but that it outputs $\widetilde{w} \in \mathfrak{S}_n$ as in equation~\eqref{tilde} satisfying (R1) and (R2).  This completes one direction of the proof.

\medskip

Now suppose $p \in \mathfrak{S}_k$ is not vexillary.  There is an occurrence $\langle 2143 \rangle$ such that
\begin{equation*}
p = \cdots \langle 2 \rangle \cdots \langle 1 \rangle (\langle 2 \rangle + 1) (\langle 2 \rangle + 2) \cdots (\langle 3 \rangle - 2)(\langle 3 \rangle - 1) \langle 4 \rangle \cdots \langle 3 \rangle \cdots.
\end{equation*}

\noindent Define $z$ to be the index such that $p(z) = \langle 1 \rangle$.  Define $w \in \mathfrak{S}_{k+1}$ by
\begin{equation*}
w(m) = \left\{ \begin{array} {c@{\quad:\quad}l}  p(m) & m \le z \text{ and } p(m) \le \langle 2 \rangle;\\
p(m) + 1 & m \le z \text{ and } p(m) > \langle 2 \rangle;\\
\langle 2 \rangle + 1 & m = z+1;\\
p(m-1) & m > z+1 \text{ and } p(m) \le \langle 2 \rangle;\\
p(m-1) + 1 & m > z+1 \text{ and } p(m) > \langle 2 \rangle. \end{array} \right.
\end{equation*}

\noindent For example, if $p = 2143$, then $w = 21354$.

If there is a reduced decomposition $\bm{j} \in R(w)$ such that $\bm{j} = \bm{j_1} \bm{i}^M \bm{j_2}$ for $\bm{i} \in R(p)$ and $M \in \mathbb{N}$, then there is a $\widetilde{w} \in \mathfrak{S}_{k+1}$ as in equation~\eqref{tilde} satisfying (R1) and (R2).  Keeping the values $\langle 1 \rangle$, $\langle 2 \rangle$, $\langle 3 \rangle$, and $\langle 4 \rangle$ as defined above, the permutation $w$ was constructed so that
\begin{equation*}
w = \cdots \langle 2 \rangle \cdots \langle 1 \rangle (\langle 2 \rangle + 1) (\langle 2 \rangle + 2) \cdots (\langle 3 \rangle - 2)(\langle 3 \rangle - 1) \langle 3 \rangle (\langle 4 \rangle + 1) \cdots (\langle 3 \rangle + 1) \cdots.
\end{equation*}

One of the values in the consecutive subsequence $(\langle 2 \rangle + 1)(\langle 2 \rangle + 2)\cdots \langle 3 \rangle$ must move to get a consecutive $p$-pattern in $\widetilde{w}$.  However, the values $\{\langle 2 \rangle, \ldots ,\langle 3 \rangle + 1\}$ appear in increasing order in $w$, and the consecutive subsequence 
\begin{equation*}
\langle 1 \rangle (\langle 2 \rangle + 1) (\langle 2 \rangle + 2) \cdots (\langle 3 \rangle - 2)(\langle 3 \rangle - 1) \langle 3 \rangle (\langle 4 \rangle + 1)
\end{equation*}
\noindent in $w$ is increasing.  Therefore, there is no way to multiply $w$ by adjacent transpositions, always eliminating an inversion, to obtain a consecutive $p$-pattern.

Hence, if $p$ is not vexillary then there exists a permutation $w$ containing a $p$-pattern such that no reduced decomposition of $w$ contains a shift of a reduced decomposition of $p$ as a factor.
\end{proof}

\begin{example}\label{vex ex}
If $w = \bm{3}14\bm{6}5\bm{2}$ and $p = 231$, with the chosen occurrence $\langle p \rangle$ in bold, the algorithm \textsf{VEX} may proceed as follows.
\begin{itemize}
\item $w_{[0]} :=  \bm{3}14\bm{6}5\bm{2}$.
\item Step~\ref{alg halt}: $x_{[0]} := 1$.
\item Step~\ref{alg small}: $w_{[0]} \mapsto w_{[0]} s_1 = 1\bm{3}4\bm{6}5\bm{2} =:w_{[1]}$.
\item Step~\ref{alg halt}: $x_{[1]} := 5$.
\item Step~\ref{alg big left}: $w_{[1]} \mapsto s_5 w_{[1]} = 1\bm{3}4\bm{5}6\bm{2} =: w_{[2]}$; $x_{[2]} := 6$.
\item Step~\ref{alg big}: $w_{[2]} \mapsto w_{[2]} s_5 = 1\bm{3}4\bm{52}6 =: w_{[3]}$.
\item Step~\ref{alg halt}: $x_{[3]} :=4$.
\item Step~\ref{alg unob right}: $w_{[3]} \mapsto w_{[3]} = 1\bm{34}5\bm{2}6 =: w_{[4]}$; $x_{[4]} := 5$.
\item Step~\ref{alg big}: $w_{[4]} \mapsto w_{[4]}s_4 = 1\bm{342}56 =: w_{[5]}$.
\item Step~\ref{alg halt}: output $1\bm{342}56$.
\end{itemize}

Therefore $\widetilde{w} = 134256 = s_5ws_1s_5s_4$, and $\widetilde{w}' = 123456$.  Keeping the notation of equation~\eqref{redw}, $\bm{h} = \emptyset$ and $M = 1$.  The unique reduced decomposition of $231$ is $12$, and indeed
\begin{equation*}
(5)\emptyset(12)^1(451) = 523451 \in R(w).
\end{equation*}
\end{example}

\begin{example}
Let $w = 21354$ and $p = 2143$.  No element of $R(w) = \{14, 41\}$ contains a shift of any element of $R(p) = \{13,31\}$ as a factor.
\end{example}

\begin{rem}
Suppose that $\bm{j} \in R(w)$ contains a shift of $\bm{i} \in R(p)$ as a factor, 
\begin{equation}\label{subword}
\bm{j} = \bm{j_1} \bm{i}^M \bm{j_2}.
\end{equation}
\noindent Then $\bm{i} \in R(p)$ can be replaced by any $\bm{i'} \in R(p)$ in equation~\eqref{subword}.
\end{rem}

Some care must be taken regarding factors in reduced decompositions.  This is clarified in the following definition and lemma, the proof of which is straightforward.

\begin{defn}
Let $w \in \mathfrak{S}_n$ and $\bm{i} \in R(w)$.  Write  $\bm{i} = \bm{a}\bm{b}\bm{c}$, where $\bm{a} \in R(u)$ and $\bm{c} \in R(v)$.  Suppose that $\bm{b}$ contains only letters in $S = \{1+M, \ldots, k-1+M\}$.  If no element of $R(u)$ has an element of $S$ as its rightmost character and no element of $R(v)$ has an element of $S$ as its leftmost character, then $\bm{b}$ is \emph{isolated} in $\bm{i}$.  Equivalently, the values $\{1+M, \ldots, k+M\}$ must appear in increasing order in $v$, and the positions $\{1+M, \ldots, k+M\}$ must comprise an increasing sequence in $u$.
\end{defn}

If $\bm{b} \in R(k \cdots 21)$ and a shift of $\bm{b}$ appears as a factor in a reduced decomposition of some permutation, then $\bm{b}$ is necessarily isolated.  This is because $\bm{b}$ has maximal reduced length in the letters $\{1, \ldots, k-1\}$, so any factor of length greater than $\binom{k}{2}$ in the letters $\{1+M, \ldots, k-1+M\}$ is not reduced.

\begin{lem}\label{isolated}
If a reduced decomposition of $w$ contains an isolated shift of a reduced decomposition of $p$, then $w$ contains the pattern $p$.
\end{lem}

The converse to Lemma~\ref{isolated} holds if $p$ is vexillary.

The characterization of vexillary in Theorem~\ref{vexthm} differs substantially from those in Proposition~\ref{vex defs}.  There is not an obvious way to prove equivalence with any of the definitions (V\ref{vex2})-(V\ref{vex6}), except via (V\ref{vex1}).  This raises the question of whether more may be understood about vexillary permutations (or perhaps other types, such as Grassmannian or dominant permutations) by studying their reduced decompositions or the permutations that contain those in question as patterns.

Theorem~\ref{vexthm} has a number of consequences, and will be used often in the subsequent sections of this paper.  Most immediately, notice that it generalizes the result of Billey, Jockusch, and Stanley mentioned earlier: $321$-avoiding permutations are exactly those whose reduced decompositions contain no long braid moves, and observe that $R(321) = \{121, 212\}$.

\section{The Commutation Relation}\label{section eq}

Recall the definition of short and long braid moves in a reduced decomposition, as well as the short and long braid relations described in equations~\eqref{shortbraid} and~\eqref{longbraid}.  It is well known that any element of $R(w)$ can be transformed into any other element of $R(w)$ by successive applications of the braid relations.

Because the short braid relation represents the commutativity of particular pairs of adjacent transpositions, the following equivalence relation is known as the \emph{commutation relation}.

\begin{defn}\label{commutation}
For a permutation $w$ and $\bm{i},\bm{j} \in R(w)$, write $\bm{i} \sim \bm{j}$ if $\bm{i}$ can be obtained from $\bm{j}$ by a sequence of short braid moves.  Let $C(w)$ be the set of commutation classes of reduced decompositions of $w$, as defined by $\sim$.
\end{defn}

\begin{example}
The commutation classes of $4231 \in \mathfrak{S}_4$ are $\{12321\}$, $\{32123\}$, and $\{13231, 31231, 13213, 31213\}$.
\end{example}

\begin{defn}\label{graph}
For a permutation $w$, the graph $G(w)$ has vertex set equal to $C(w)$, and two vertices share an edge if there exist representatives of the two classes that differ by a long braid move.
\end{defn}

Elnitsky gives a very elegant representation of this graph in \cite{elnitsky}, which will be discussed in depth in Section~\ref{section polygon}.  A consequence of his description, although not difficult to prove independent of his work, is the following.

\begin{prop}\label{graph facts}
The graph $G(w)$ is connected and bipartite.
\end{prop}

\begin{proof}
See \cite{elnitsky}.
\end{proof}

Despite Proposition~\ref{graph facts}, much remains to be understood about the graph $G(w)$.  For example, even the size of the graph for $w_0$ (that is, the number of commutation classes for the longest element) is unknown.

Billey, Jockusch, and Stanley characterize all permutations with a single commutation class, and hence whose graphs are a single vertex, as $321$-avoiding permutations.  A logical question to ask next is: for what permutations does each reduced decomposition contain at most one long braid move?  More restrictively: what if this long braid move is required to be a specific shift of $121$ or $212$?  Moreover, what are the graphs in these cases?

\begin{defn}
Let $U_n = \{w \in \mathfrak{S}_n : \text{no } \bm{j} \in R(w) \text{ has two long braid moves}\}$.
\end{defn}

\begin{thm}\label{1 lbm}
$U_n$ is the set of permutations such that every $321$-pattern in $w$ has the same maximal element and the same minimal element.
\end{thm}

\begin{proof}
Assume $w$ has a $321$-pattern.  Suppose that every occurrence of $321$ in $w$ has $\langle 3 \rangle = x$ and $\langle 1 \rangle = y$.  Suppose that $\bm{j} \in R(w)$ has at least one long braid move.  Choose $k$ so that $j_k j_{k+1} j_{k+2}$ is the first such.  Each adjacent transposition in a reduced decomposition increases the length of the product.  Then by the supposition,
\begin{equation*}
s_{j_{k+2}} s_{j_{k+1}} s_{j_k} \cdots s_{j_1} w
\end{equation*}
\noindent is $321$-avoiding, so $j_{k+3} \cdots j_\ell$ has no long braid moves.  It remains only to consider when $j_{k+2} j_{k+3} j_{k+4}$ is also a long braid move.  The only possible reduced configurations for such a factor $j_k j_{k+1} j_{k+2} j_{k+3} j_{k+4}$ are shifts of $21232$ and $23212$.  If either of these is not isolated, then it is part of a shift of $212321$, $321232$, $123212$, or $232123$.  Notice that  

\begin{itemize}
\item $212321, 321232, 123212, 232123 \in R(4321)$;
\item $23212 \in R(4312)$;
\item $21232 \in R(3421)$.
\end{itemize}

\noindent If $j_k j_{k+1} j_{k+2} j_{k+3} j_{k+4}$ is isolated in $\bm{j}$ then $w$ contains a $4312$- or $3421$-pattern by Lemma~\ref{isolated}.  Otherwise, $w$ contains a $4321$-pattern.  However, every $321$-pattern in $w$ has $\langle 3 \rangle = x$ and $\langle 1 \rangle = y$.  Therefore $j_kj_{k+1}j_{k+2}$ is the only long braid move in $\bm{j}$, so $w \in U_n$.

\smallskip

Now let $w$ be an element of $U_n$.  If $w$ has two $321$-patterns that do not have the same maximal element and the same minimal element, then they intersect at most once or they create a $4321$-, $4312$-, or $3421$-pattern.  These three patterns are vexillary.  Thus by Theorem~\ref{vexthm} and the examples above, containing one of these patterns would imply that some element of $R(w)$ has more than one long braid move.  If the two $321$-patterns intersect at most once, their union may be a non-vexillary pattern, so Theorem~\ref{vexthm} does not necessarily apply.  However, a case analysis shows that it is possible to shorten $w$ by adjacent transpositions and make one $321$-pattern increasing (via a long braid move) without destroying the other $321$-pattern.  Thus an element of $R(w)$ would have more than one long braid move, contradicting $w \in U_n$.
\end{proof}

\begin{defn}
Let $U_n(j)$ consist of permutations with some $321$-pattern, where every long braid move that occurs must be $j(j+1)j$ or $(j+1)j(j+1)$.
\end{defn}

\begin{cor}\label{1 ulbm}
$U_n(j) = \{w \in \mathfrak{S}_n : w$ has a unique $321$-pattern and $\langle 2 \rangle = j+1\}$.  If $w$ has a unique $321$-pattern, then $w(\langle 2 \rangle) = \langle 2 \rangle$.
\end{cor}

\begin{proof}
A unique $321$-pattern implies that $\{1, \ldots, \langle 2 \rangle - 1\} \setminus \langle 1 \rangle$ all appear to the left of $\langle 2 \rangle$ in $w(1) \cdots w(n)$, and $\{\langle 2 \rangle + 1, \ldots, n\} \setminus \langle 3 \rangle$ all appear to the right of $\langle 2 \rangle$, so the second statement follows.

Consider the long braid moves that may appear for elements of $U_n \supset U_n(j)$.  Let $w \in U_n$ have $k$ distinct $321$-patterns.  By Theorem~\ref{1 lbm}, these form a pattern $p = (k+2)23\cdots k(k+1)1 \in \mathfrak{S}_{k+2}$ in $w$.  The permutation $p$ is vexillary, so there exists $M \in \mathbb{N}$ and a reduced decomposition $\bm{j_1} \bm{i}^M \bm{j_2} \in R(w)$ for each $\bm{i} \in R(p)$.  There are elements in $R(p)$ with long braid moves $i(i+1)i$ for each $i \in [1,k]$.  For example, $12 \cdots k(k+1)k \cdots 21 \in R(p)$.  Therefore, if $w \in U_n(j)$, then $k = 1$, so $w$ has a unique $321$-pattern.

Suppose that $w$ has a unique $321$-pattern.  Because $w(\langle 2 \rangle) = \langle 2 \rangle$, the only possible long braid moves in reduced decompositions of $w$ are $(\langle 2 \rangle -1)\langle 2 \rangle (\langle 2 \rangle - 1)$ or $\langle 2 \rangle (\langle 2 \rangle - 1) \langle 2 \rangle$.
\end{proof}

\begin{cor}\label{U_n graph}
If $w \in U_n$ and $w$ has $k$ distinct $321$-patterns, then $|C(w)| = k+1$ and the graph $G(w)$ is a path of $k+1$ vertices connected by $k$ edges.
\end{cor}

\begin{proof}
Because $w$ contains the pattern $p = (k+2)23\cdots k(k+1)1 \in \mathfrak{S}_{k+2}$, there is a subgraph of $G(w)$ that is a path of $k+1$ vertices connected by $k$ edges.  Since $p$ accounts for all of the $321$-patterns in $w$, this is all of $G(w)$.
\end{proof}

\begin{cor}\label{U_n(j) graph}
If $w \in U_n(j)$, then $|C(w)| = 2$ and the graph $G(w)$ is a pair of vertices connected by an edge.
\end{cor}

\section{Elnitsky's Polygon}\label{section polygon}

In his doctoral thesis and in \cite{elnitsky}, Elnitsky developed a bijection between commutation classes of reduced decompositions of $w \in \mathfrak{S}_n$ and rhombic tilings of a particular $2n$-gon $X(w)$.  This bijection leads to a number of interesting questions about tilings of $X(w)$ and their relations to the permutation $w$ itself.  A number of these ideas are studied in this and the following section.

\begin{defn}\label{polygon}
For $w \in \mathfrak{S}_n$, let $X(w)$ be the $2n$-gon with all sides of unit length such that
\begin{enumerate}
\item Sides of $X(w)$ are labeled $1, \ldots, n, w(n), \ldots, w(1)$ in order;
\item The portion labeled $1, \ldots, n$ is convex; and
\item Sides with the same label are parallel.
\end{enumerate}
Orient the polygon so that the edge labeled $1$ lies to the left of the top vertex and the edge labeled $w(1)$ lies to its right.  This is \emph{Elnitsky's polygon}.
\end{defn}

\begin{figure}[htbp]
\centering
\epsfig{file=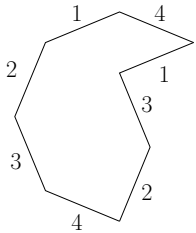, width=1in}
\caption{The polygon $X(4132)$.}
\end{figure}

\begin{example}
For $w_0 \in \mathfrak{S}_n$, the polygon $X(w_0)$ is a centrally symmetric $2n$-gon.
\end{example}

\begin{defn}
The hexagon $X(321)$ can be tiled by rhombi with sides of unit length in exactly two ways.  Each of these tilings is the \emph{flip} of the other.
\end{defn}

\begin{defn}
Let $T(w)$ be the set of tilings of $X(w)$ by rhombi with sides of unit length.  Define a graph $G'(w)$ with vertex set $T(w)$, and connect two tilings by an edge if they differ by a flip of the tiling of a single sub-hexagon.
\end{defn}

Unless otherwise indicated, the term \emph{tiling} refers to an element of $T(w)$.

\begin{thm}[Elnitsky]\label{elthm}
The graphs $G(w)$ and $G'(w)$ are isomorphic.
\end{thm}

Henceforth, both graphs will be denoted $G(w)$.

Before discussing new results related to this polygon, it is important to understand Elnitsky's bijection, outlined in the following algorithm.  A more thorough treatment appears in \cite{elnitsky}.

\medskip

\begin{myquote}

\setlength{\leftmargini}{0in}
\setlength{\leftmarginii}{.18in}
\setlength{\leftmarginiii}{.18in}
\setlength{\leftmarginiv}{.18in}

\noindent \textsf{\textbf{Algorithm} ELN}

\nopagebreak[4]

\noindent \textsf{INPUT: $T \in T(w)$.}

\noindent \textsf{OUTPUT: An element of $C_T \in C(w)$.}

\smallskip

\begin{enumerate}
\setcounter{enumi}{-1}
\renewcommand{\labelenumi}{\textsf{Step \arabic{enumi}}.}

\item Set the polygon $P_{[0]} := X(w)$, the string $\bm{j_{[0]}} := \emptyset$, and $i:=0$.
\item\label{elnstart} If $P_{[i]}$ has no area, then \textsf{OUTPUT $\bm{j_{[i]}}$}.
\item\label{elnchoice} There is at least one tile $t_i$ that shares two edges with the right side of $P_{[i]}$.  
\item If $t_i$ includes the $j^{\text{th}}$ and $(j+1)^{\text{st}}$ edges from the top along the right side of $P_{[i]}$, set $\bm{j_{[i+1]}} := j\bm{j_{[i]}}$.
\item Let $P_{[i+1]}$ be $P_{[i]}$ with the tile $t_i$ removed.  Set $i:=i+1$ and \textsf{GOTO Step~\ref{elnstart}}.
\end{enumerate}
\end{myquote}

\textsf{ELN} yields the entire commutation class because of the choice of tile in Step~\ref{elnchoice}.

\medskip

\begin{figure}[htbp]
\centering
\epsfig{file=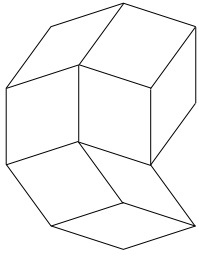, width=1in}
\caption{A tiling in $T(53241)$.}\label{tilingfig}
\end{figure}

\begin{example}
The tiling in Figure~\ref{tilingfig} corresponds to the equivalence class consisting solely of the reduced decomposition $12343212 \in R(53241)$.
\end{example}

\begin{cor}\label{subgraph}
If $p$ is vexillary and $w$ contains a $p$-pattern, then $G(p)$ is a subgraph of $G(w)$.
\end{cor}

\begin{proof}
This follows from Theorems~\ref{vexthm} and~\ref{elthm}.
\end{proof}

Elnitsky's correspondence, described in \textsf{ELN}, combined with Theorem~\ref{1 lbm} and Corollary~\ref{1 ulbm}, indicates that any tiling of $X(w)$ for $w \in U_n$ has at most one sub-hexagon (every tiling has exactly one sub-hexagon if $w$ is not $321$-avoiding).  Moreover, the sub-hexagon has the same vertical position for all elements of $U_n(j)$.

Under certain circumstances, the polygon $X(w)$ for $w \in \mathfrak{S}_n$ can be rotated or reflected to give a polygon $X(w')$ for another $w' \in \mathfrak{S}_n$.

\begin{cor}\label{reflect}
Let $w = w(1) \cdots w(n)$ and $w^R = w(n) \cdots w(1)$.  Then $|C(w)| = |C(w^R)|$ and $G(w) \simeq G(w^R)$.
\end{cor}

\begin{cor}\label{rotate}
Let $w = w(1) \cdots w(n)$.  If $w(1) = n, w(2) = n-1, \ldots, w(i) = n+1 - i$, then $|C(w)| = |C(w^{(i)})|$and $G(w) \simeq G(w^{(i)})$ where
\[w^{(i)} = (w(i+1) + i)(w(i+2)+i) \cdots (w(n)+i)i(i-1) \cdots 21\]
\noindent and all entries are modulo $n$.  Likewise, if $w(n) = 1, w(n-1) = 2, \ldots, w(n-j+1) = j$, then $|C(w)| = |C(w_{(j)})|$ and $G(w) \simeq G(w_{(j)})$ where
\[w_{(j)} = n(n-1) \cdots (n-j+1) (w(1)-j) (w(2)-j) \cdots (w(n-j)-j)\]
\noindent and all entries are modulo $n$.
\end{cor}

Elnitsky's result interprets the commutation classes of $R(w)$ as rhombic 
tilings of $X(w)$, with long braid moves represented by flipping sub-hexagons.  
The following theorem utilizes this interpretation, and demonstrates that the 
number of commutation classes of a permutation is monotonically increasing 
with respect to pattern containment, thus generalizing one aspect of 
Corollary~\ref{subgraph}.  Note that $p$ is not required to be vexillary in 
Theorem~\ref{monotone}, unlike in Theorem~\ref{vexthm}.

\begin{thm}\label{monotone}
If $w$ contains the pattern $p$, then $|C(w)| \ge |C(p)|$.
\end{thm}

\begin{proof}
Consider a tiling $T \in T(p)$.  This represents a commutation class of $R(p)$.  For an ordering of the tiles in $T$ as defined by \textsf{ELN}, label the tile $t_0$ by $\ell(p)$, the tile $t_1$ by $\ell(p) - 1$, and so on.  If the tile with label $r$ corresponds to the adjacent transposition $s_{i_r}$, then $i_1 \cdots i_{\ell(p)} \in R(p)$.

\medskip

\begin{myquote}

\setlength{\leftmargini}{0in}
\setlength{\leftmarginii}{.18in}
\setlength{\leftmarginiii}{.18in}
\setlength{\leftmarginiv}{.18in}

\noindent \textsf{\textbf{Algorithm} MONO}

\nopagebreak[4]

\noindent \textsf{INPUT: $w$ containing the pattern $p$ and $T \in T(p)$ with tiles labeled as described.}

\noindent \textsf{OUTPUT: $T' \in T(w)$.}

\smallskip

\begin{enumerate}
\setcounter{enumi}{-1}
\renewcommand{\labelenumi}{\textsf{Step \arabic{enumi}}.}

\item Set $w_{[0]}:=w$, $p_{[0]}:=p$, $T_{[0]}:=T$, $T'_{[0]}:=\emptyset$, and $i:=0$.
\item\label{monostart} If $p_{[i]}$ is the identity permutation, then define $T'_{[i+1]}$ to be the tiles of $T'_{[i]}$ together with any tiling of $X(w_{[i]})$.  \textsf{OUTPUT $T'_{[i+1]}$}.
\item Let $j_{[i]}$ be such that the tile labeled $\ell(p)-i$ includes edges $p_{[i]}(j_{[i]})$ and $p_{[i]}(j_{[i]}+1)$.  Note that $p_{[i]}(j_{[i]}) > p_{[i]}(j_{[i]}+1)$.
\item Define $r<s$ so that $w_{[i]}(r) = \langle p_{[i]}(j_{[i]}) \rangle$ and $w_{[i]}(s) = \langle p_{[i]}(j_{[i]}+1) \rangle$.  Note that $w_{[i]}(t)$ is a non-pattern entry for $t \in (r,s)$.
\item Let $v_{[i]}$ be the permutation defined by
\begin{equation*}
w_{[i]}(t) \mapsto \left\{ \begin{array} {c@{\quad:\quad}l}
w_{[i]}(t) & t<r \text{ or } t>s\\
\widetilde{w}_{[i]}(t) & r \le t \le s
\end{array}
\right.
\end{equation*}
\noindent where $(\widetilde{w}_{[i]}(r), \ldots, \widetilde{w}_{[i]}(s))$ is $\{w_{[i]}(r), \ldots, w_{[i]}(s)\}$ in increasing order.
\item Set $w_{[i+1]} := w_{[i]}v_{[i]}$, and notice that $\ell(w_{[i+1]}) = \ell(w_{[i]}) - \ell(v_{[i]})$.
\item\label{monochoice} The right boundaries of $X(w_{[i+1]})$ and $X(w_{[i]})$ differ only in the $r^{\text{th}}, \ldots, s^{\text{th}}$ edges, and the left side of this difference (part of the boundary of $X(w_{[i+1]})$) is convex.  Therefore, this difference has a rhombic tiling $t_{[i]}$.  Define $T'_{[i+1]}$ to be the tiles in $T'_{[i]}$ together with the tiles in $t_{[i]}$.
\item Set $i:=i+1$ and \textsf{GOTO Step~\ref{monostart}}.
\end{enumerate}
\end{myquote}

\medskip

The algorithm \textsf{MONO} takes a tiling $T \in T(p)$ and outputs one of possibly several tilings $T' \in T(w)$ due to the choice in Steps~\ref{monostart} and~\ref{monochoice}.  A tiling $T' \in T(w)$ so obtained can only come from this $T$, although possibly with more than one labeling of the tiles.  However, this labeling of the tiles merely reflects the choice of a representative from the commutation class, so indeed $|T(w)| \ge |T(p)|$, and $|C(w)| \ge |C(p)|$.
\end{proof}

\begin{figure}[htbp]
\centering
\parbox{.8in}{\epsfig{file=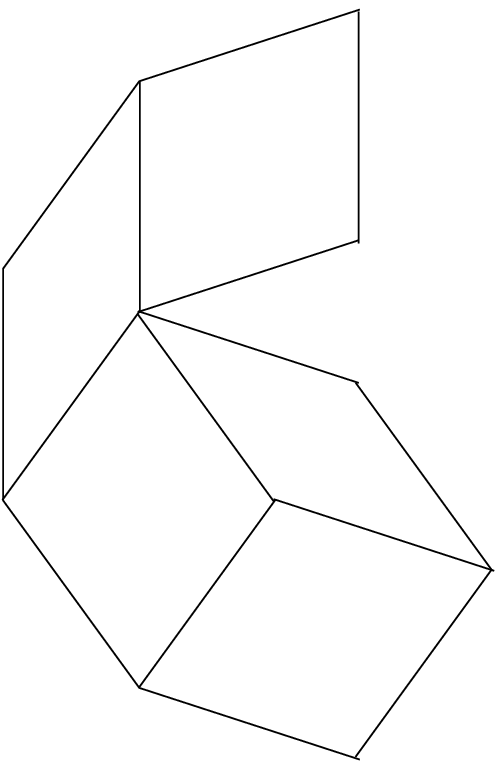, width=.8in}} \hspace{.25in} $\Longrightarrow$ \hspace{.25in} \parbox{1.04in}{\epsfig{file=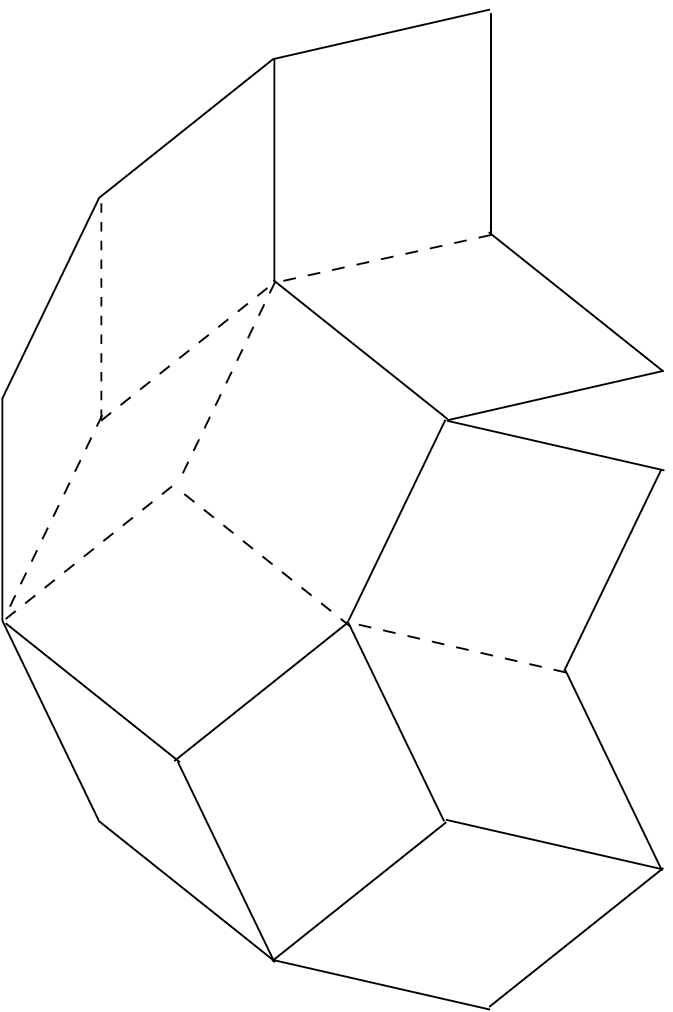, width=1.04in}}\\
\parbox{.8in}{\epsfig{file=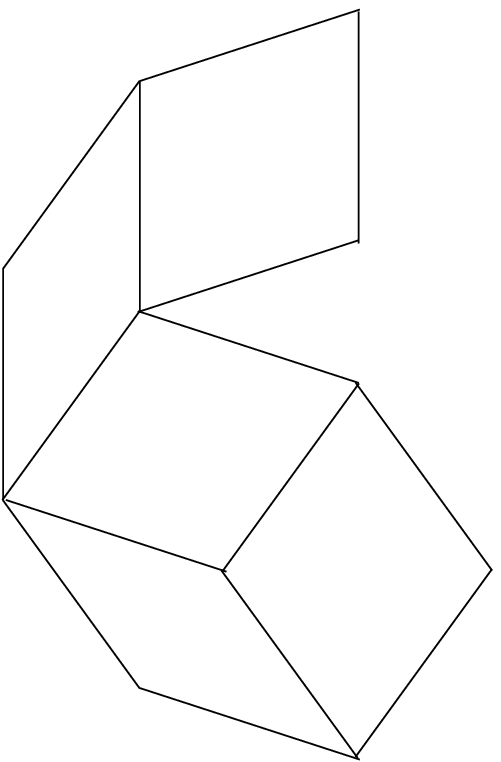, width=.8in}} \hspace{.25in} $\Longrightarrow$ \hspace{.25in} \parbox{1.04in}{\epsfig{file=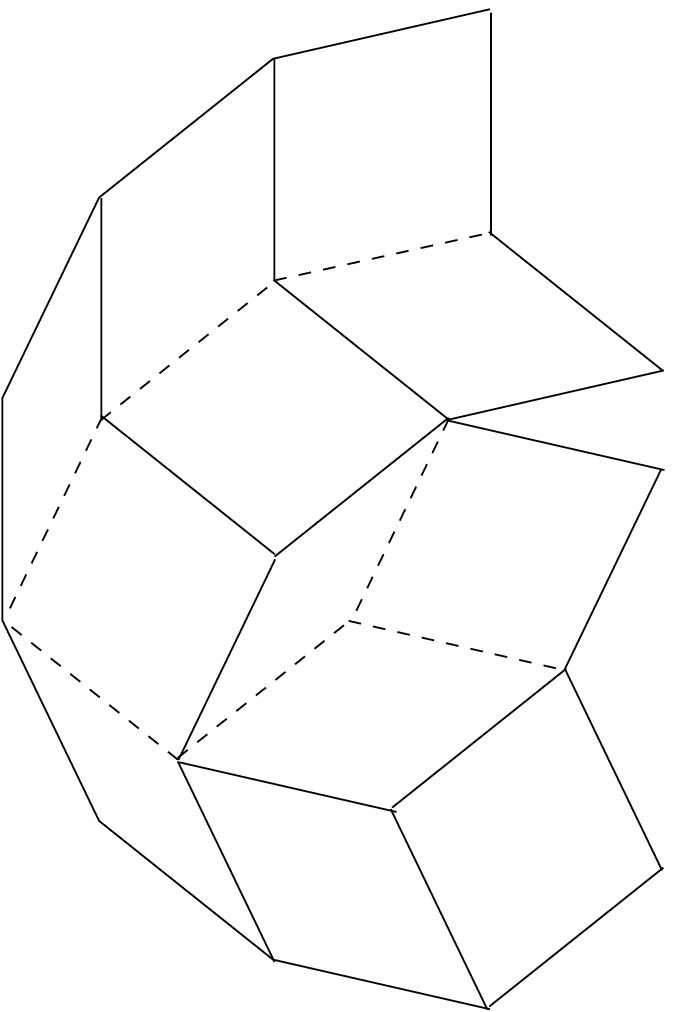, width=1.04in}}
\caption{The output of \textsf{MONO}, given each of the two elements of $T(31542)$.  The dotted lines indicate the choice of tiling in Steps~\ref{monostart} and~\ref{monochoice} of the algorithm.}\label{mono ex}
\end{figure}

\begin{example}
Let $p = 31542$ and $w = 4617352$.  The pattern $p$ occurs in $w$ as $\langle p \rangle = 41752$.  Figure~\ref{mono ex} depicts the output of \textsf{MONO}, given the two tilings of $X(p)$.  
\end{example}

\section{The Poset of Tilings}\label{section poset}

Elnitsky's bijection considers the rhombic tilings of the polygon $X(w)$.  Rhombi are a special case of a more general class of objects known as zonotopes.

\begin{defn}
A polytope is a \emph{$d$-zonotope} if it is the projection of a regular $n$-cube onto a $d$-dimensional subspace.
\end{defn}

Centrally symmetric convex polygons are exactly the $2$-zonotopes.  These necessarily have an even number of sides.

\begin{defn}
A \emph{zonotopal tiling} of a polygon is a tiling by centrally symmetric convex polygons.
\end{defn}

\begin{defn}
Let $Z(w)$ be the set of zonotopal tilings of Elnitsky's polygon.  Rhombi are centrally symmetric, so $T(w) \subseteq Z(w)$.
\end{defn}

\begin{figure}[htbp]
\centering
\epsfig{file=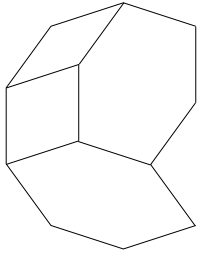, width=1in}
\caption{A tiling in $Z(53241)$.}\label{ztilingfig}
\end{figure}

\begin{thm}\label{2k-gon}
There is a tiling in $Z(w)$ containing a $2k$-gon with sides parallel to the sides labeled $i_1 < \cdots < i_k$ if and only if $i_k \cdots i_1$ is an occurrence of $k \cdots 1$ in $w$.
\end{thm}

\begin{proof}
Since the tiles are convex, a $2k$-gon in the tiling with sides as described has right side labeled $i_k, \ldots, i_1$ from top to bottom and left side labeled $i_1, \ldots, i_k$ from top to bottom.  Therefore Elnitsky's bijection shows that this tile (or rather, any decomposition of it into rhombi) transforms the sequence $(i_1, \ldots, i_k)$ into $(i_k, \ldots, i_1)$.  Reduced decompositions have minimal length, so no inversions can be ``undone'' by subsequent adjacent transpositions.  Therefore $i_k \cdots i_1$ must be an occurrence of $k \cdots 1$ in $w$.

Conversely, suppose that $i_k \cdots i_1$ is an occurrence of the vexillary pattern $k \cdots 1$ in $w$.   For a decreasing pattern, the algorithm \textsf{VEX} can be modified slightly to produce $\widetilde{w}$ as in equation~\eqref{tilde}, where the consecutive occurrence $\langle k \cdots 1 \rangle$ is $i_k \cdots i_1$.  Let $\bm{i} \in R(k \cdots 1)$ and $(I_q \cdots I_1)\bm{h}\bm{i}^M(J_r \cdots J_1) \in R(w)$ for $\bm{h} \in R(\widetilde{w}')$.  Removing the rhombi that correspond to $s_{J_r} \cdots s_{J_1}$ yields the polygon $X(\widetilde{w})$, and the rhombi that correspond to $\bm{i}^M$ form a sub-$2k$-gon with sides parallel to the sides labeled $\{i_1, \ldots, i_k\}$ in $X(w)$.
\end{proof}

Less specifically, Theorem~\ref{2k-gon} states that a tiling in $Z(w)$ can contain a $2k$-gon if and only if $w$ has a decreasing subsequence of length $k$.


Using a group theoretic argument, Pasechnik and Shapiro showed in 
\cite{pasechnik} that no element of $Z(n \cdots 21)$ consists entirely of 
hexagons for $n>3$.  Their result states that at least one rhombus must be 
present in a hexagonal/rhombic tiling.  Kelly and Rottenberg had previously 
obtained a better bound in terms of arrangements of pseudolines in \cite{kelly}.

Working with reduced decompositions and Elnitsky's polygons yields a different proof that no element of $Z(n \cdots 21)$ can consist of entirely hexagonal tiles for $n>3$, and generalizes the result to other types of tiles.  Theorem~\ref{2k tiles} is, in a sense, a counterpart to Theorem~\ref{2k-gon}.

\begin{thm}\label{2k tiles}
Let $w_0$ be the longest element in $\mathfrak{S}_n$.  There is a tiling $Z \in Z(w_0)$ consisting entirely of $2k$-gons if and only if one of the following is true:
\begin{enumerate}
\item $k = 2$; or
\item $k = n$.
\end{enumerate}
\end{thm}

\begin{proof}
If $k=2$, the result holds for all $n$: every $X(w_0)$ can be tiled by rhombi with unit side length.  For the remainder of the proof, assume that $k \ge 3$.

Suppose that there exists $Z \in Z(w_0)$ consisting entirely of $2k$-gons.  Then there exists $\bm{j}^{M_1} \cdots \bm{j}^{M_r} \in R(w_0)$, for some $\bm{j} \in R(k \cdots 21)$.  Because $k \cdots 21$ is the longest element in $\mathfrak{S}_k$, the factor $\bm{j}^{M_i}$ can have any of $\{1+M_i, \ldots, k-1 + M_i\}$ at either end.  Thus, if $M_i = M_j = M$ for $i<j$, then $M \pm (k-1) \in \{M_{i+1}, \ldots, M_{j-1}\}$.  Because $M_i \in [0,n-k]$, there is at most one shift by $0$ and at most one shift by $n-k$.  In fact, there is exactly one of each of these shifts, since all of $X(w_0)$ must be tiled and the shifts correspond to the vertical placement of the $2k$-gons.

Consider the tile corresponding to $\bm{j}^0$.  This $2k$-gon is as high vertically as possible and it is the only tile placed so high.  Thus it shares the top vertex and its incident sides with $X(w_0)$.  That is, two of the tile's sides are labeled $1$ and $n$.  Similarly, the tile for $\bm{j}^{n-k}$ also has two sides labeled $1$ and $n$.

By Theorem~\ref{2k-gon}, each of these tiles corresponds to a $k \cdots 21$-pattern with $\langle k \rangle = n$ and $\langle 1 \rangle = 1$.  However, once the value $1$ is to the right of the value $n$, it remains to the right as adjacent positions are transposed to lengthen the permutation.  Therefore the two patterns, the tiles, and their shifts must be equal: $0 = n-k$.

Indeed, there is always a tiling $Z \in Z(w_0)$ consisting of a single $2n$-gon.
\end{proof}


There is a poset $P(w)$ that arises naturally when studying $Z(w)$.

\begin{defn}
For a permutation $w$, let the poset $P(w)$ have elements equal to the zonotopal tilings $Z(w)$, partially ordered by reverse edge inclusion.  
\end{defn}

\begin{example}
In the poset $P(53241)$, the tiling in Figure~\ref{tilingfig} is less than the tiling in Figure~\ref{ztilingfig}.
\end{example}

\begin{rem}\label{P(w_0) max}
For the longest element $w_0 \in \mathfrak{S}_n$, the poset $P(w_0)$ has a maximal element equal to the tiling in $Z(w_0)$ that consists of a single $2n$-gon.
\end{rem}

\begin{rem}\label{posetgraph}
The minimal elements of $P(w)$ are the rhombic tilings, which are the vertices of the graph $G(w)$.  Moreover, edges in the graph $G(w)$ correspond to flipping a single sub-hexagon in the tiling.  Therefore these edges correspond to the elements of $P(w)$ that cover the minimal elements.
\end{rem}

The relationships in Remark~\ref{posetgraph} are immediately apparent.  Another relationship is not as obvious.  This follows from a result of Boris Shapiro, Michael Shapiro, and Alek Vainshtein in \cite{ssv}.

\begin{lem}[Shapiro-Shapiro-Vainshtein]\label{ssvlem}
The set of all $4$- and $8$-cycles in $G(w)$ form a system of generators for the first homology group $H_1(G(w), \mathbb{Z}/2\mathbb{Z})$.
\end{lem}

Additionally, Anders Bj\"{o}rner noted that gluing $2$-cells into those $4$- and $8$-cycles yields a simply connected complex (\cite{bjorner}).

In \cite{ssv}, Lemma~\ref{ssvlem} is stated only for $w = w_0$.  However, the proof easily generalizes to all $w \in \mathfrak{S}_n$.  A straightforward argument demonstrates that a $4$-cycle in $G(w)$ corresponds to $Z \in Z(w)$ with rhombi and two hexagons, and an $8$-cycle corresponds to $Z \in Z(w)$ with rhombi and an octagon.  These are exactly the elements of $P(w)$ which cover those that correspond to edges of $G(w)$.

\begin{cor}\label{level2}
The elements of $P(w)$ that cover the elements (corresponding to edges of $G(w)$) covering the minimal elements (corresponding to vertices of $G(w)$) correspond to a system of generators for the first homology group $H_1(G(w), \mathbb{Z}/2\mathbb{Z})$.
\end{cor}

Little is known about the structure of the graph $G(w)$ for arbitrary $w$.  However, in some cases a description can be given via Theorem~\ref{2k-gon} and Lemma~\ref{ssvlem}.

\begin{thm}\label{chainthm}
The following statements are equivalent for a permutation $w$:

\begin{enumerate}
\item\label{treecon} $G(w)$ is a tree;
\item\label{chaincon} $G(w)$ is a path (that is, no vertex has more than two incident edges);
\item\label{posetcon} The maximal elements of $P(w)$ cover the minimal elements;
\item\label{patterncon} $w$ is $4321$-avoiding and any two $321$-patterns intersect at least twice.
\end{enumerate}
\end{thm}

\begin{proof}
(\ref{treecon}) $\Leftrightarrow$ (\ref{posetcon}) by Corollary~\ref{level2}.  From Theorem~\ref{2k-gon} and the discussion preceding Corollary~\ref{level2}, an $8$-cycle in the graph is equivalent to having a $4321$-pattern.  Similarly, a $4$-cycle is equivalent to two sub-hexagons whose intersection has zero area, so some reduced decomposition has two disjoint long braid moves.  This implies that two $321$-patterns intersect in at most one position.  Therefore (\ref{treecon}) $\Leftrightarrow$ (\ref{patterncon}).

Finally, suppose that $G(w)$ is a tree and a vertex has three incident edges.  The corresponding tiling has at least three sub-hexagons.  However, it is impossible for every pair of these to overlap.  This contradicts (\ref{treecon}) $\Leftrightarrow$ (\ref{posetcon}) $\Leftrightarrow$ (\ref{patterncon}), so (\ref{treecon}) $\Leftrightarrow$ (\ref{chaincon}).
\end{proof}

If $C_n$ is the set of all $w \in \mathfrak{S}_n$ for which $G(w)$ is a path, then $U_n(j) \subseteq U_n \subseteq C_n$ by Corollaries~\ref{U_n graph} and~\ref{U_n(j) graph}.

Following convention, the unique maximal element in a poset, if it exists, is denoted $\widehat{1}$.  Remark~\ref{P(w_0) max} noted that the poset $P(w_0)$ has a $\widehat{1}$.  In fact, there are other $w$ for which $P(w)$ has a $\widehat{1}$, as described below.

\begin{thm}\label{max elt}
The poset $P(w)$ has a $\widehat{1}$ if and only if $w$ is $4231$-, $4312$-, and $3421$-avoiding.
\end{thm}

\begin{proof}
The definition of the poset $P(w)$ and Theorem~\ref{2k-gon} indicate that $P(w)$ has a $\widehat{1}$ if and only if the union of any two decreasing subsequences that intersect at least twice is itself a decreasing subsequence.

Suppose there are decreasing subsequences in $w$ of lengths $k_1,k_2 \ge 3$ that intersect $i \ge 2$ times, for $i < k_1,k_2$.  Let $k = i+1$, and choose a $k+1$ element subsequence of $\langle k_1 \cdots 1 \rangle \cup \langle k_2 \cdots 1 \rangle$ that includes $\langle k_1 \cdots 1 \rangle \cap \langle k_2 \cdots 1 \rangle$ and one more element from each descending subsequence.  Let $p \in \mathfrak{S}_{k+1}$ be the resulting pattern.  No $\widehat{1}$ in $P(w)$ is equivalent to there being subsequences so that
\begin{equation*}
p = (k+1)k \cdots (j+2)j(j+1)(j-1) \cdots 21
\end{equation*}
\noindent for some $j \in [1, k]$.  

There are two ways to place a $2k$-gon in a zonotopal tiling of $X(p)$, but these overlapping $2k$-gons do not both lie in any larger centrally symmetric polygon.  The permutation $p$ is always vexillary, so Theorem~\ref{vexthm} implies that $P(w)$ will not have a $\widehat{1}$ if $w$ contains such a $p$.

Therefore, considering the permutation $p$ for each possible $j$, the poset $P(w)$ has a $\widehat{1}$ if and only if $w$ is $4231$-, $4312$-, and $3421$-avoiding.
\end{proof}

The permutations for which $P(w)$ has a $\widehat{1}$ have recently been 
enumerated by Toufik Mansour in \cite{mansourpreprint}.

\section{The Freely Braided Case}\label{section fb}

Although the graph $G(w)$ and poset $P(w)$ are not known in general, there is a class of permutations for which these objects can be completely described.  This paper concludes with a study of this special case.

In \cite{green1} and \cite{green2}, Green and Losonczy introduce and study ``freely braided'' elements in simply laced Coxeter groups.  In the case of type $A$, these are as follows.

\begin{defn}
A permutation $w$ is \emph{freely braided} if every pair of distinct $321$-patterns in $w$ intersects at most once.
\end{defn}

Equivalently, $w$ is freely braided if and only if $w$ is $4321$-, $4231$-, $4312$-, and $3421$-avoiding.  The poset of a freely braided permutation has a unique maximal element by Theorem~\ref{max elt}.

\begin{example}
The permutation $35214$ is not freely braided because $321$ and $521$ are both occurrences of the pattern $321$, and they intersect twice.  The permutation $52143$ is freely braided.
\end{example}

Mansour enumerates freely braided permutations in \cite{mansour}.

In \cite{green1}, Green and Losonczy show that a freely braided $w$ with $k$ distinct $321$-patterns has
\begin{equation}\label{fb size}
|C(w)| = 2^k.
\end{equation}
\noindent Moreover, in \cite{green2} they show the following fact for any simply laced Coxeter group, here stated only for type $A$.

\begin{prop}[Green-Losonczy]
If a permutation $w$ is freely braided with $k$ distinct $321$-patterns, then there exists $\bm{i} \in R(w)$ with $k$ disjoint long braid moves.
\end{prop}

\begin{rem}\label{fb tiles}
This means that there is a tiling of $X(w)$ with $k$ sub-hexagons, none of which overlap.  Furthermore, equation~\eqref{fb size} implies that flipping any sequence of these sub-hexagons does not yield any new sub-hexagons.  Hence \emph{every} tiling of $X(w)$ has exactly $k$ sub-hexagons, none of which overlap, and $\widehat{1}$ in $P(w)$ corresponds to the zonotopal tiling with rhombi and $k$ hexagons.
\end{rem}

From Remark~\ref{fb tiles}, the structures of the graph $G(w)$ and the poset $P(w)$ are clear for a freely braided permutation $w \in \mathfrak{S}_n$.

\begin{thm}
Let $w$ be freely braided with $k$ distinct $321$-patterns.  The graph $G(w)$ is the graph of the $k$-cube, and the poset $P(w)$ is isomorphic to the face lattice of the $k$-cube without its minimal element.
\end{thm}

\begin{figure}[htbp]
\centering
\epsfig{file=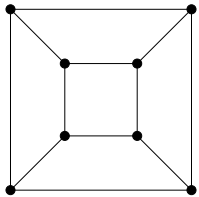, width=1in}
\caption{$G(243196587)$.}\label{3cubefig}
\end{figure}

\begin{figure}[htbp]
\centering
\epsfig{file=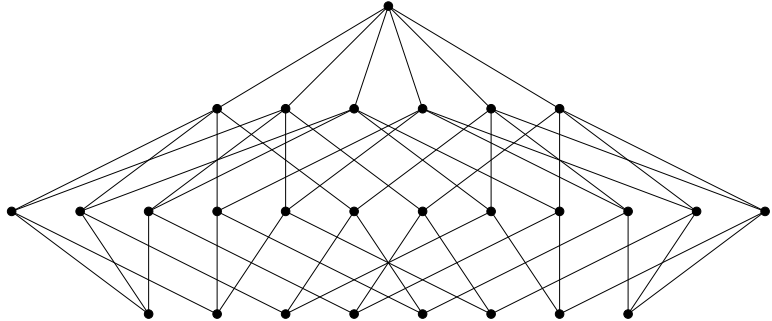, width=4in}
\caption{$P(243196587)$.}\label{3posetfig}
\end{figure}

\begin{example}
The permutation $243196587$ is freely braided.  Its three $321$-patterns are $431$, $965$, and $987$.  Figures~\ref{3cubefig} and~\ref{3posetfig} depict its graph and poset.
\end{example}

\section{Acknowledgments}\label{section ack}

Particular thanks are due to Richard Stanley for his continued guidance and for the suggestion to study reduced decompositions.  Anders Bj\"{o}rner provided helpful advice and discussion, and together with Richard Stanley coordinated the semester on algebraic combinatorics at the Institut Mittag-Leffler, during which much of the research for this paper occurred.  Thanks are also owed to John Stembridge for his referral to the work of Elnitsky and Green and Losonczy.  Finally, the thoughtful suggestions of two anonymous referees have been greatly appreciated.

\end{document}